\documentclass[a4paper,11pt]{article}
\usepackage[T1]{fontenc}
\usepackage[utf8]{inputenc}
\usepackage{amssymb}
\usepackage{amsmath}
\usepackage{graphicx}
\usepackage{geometry}
\usepackage{color}
\newcommand{\highlight}[1]{#1}
\newtheorem{Theorem}{Theorem}
\newtheorem{Lemma}[Theorem]{Lemma}
\newtheorem{Definition}[Theorem]{Definition}
\newtheorem{Example}[Theorem]{Example}
\newenvironment{Proof}
{\begin{trivlist}\item[]{{\sc Proof.}}}{\hfill{$\square$}\noindent\end{trivlist}
}

\begin{document}

\title{Generalized roll-call model for the Shapley-Shubik index}

\author{Sascha Kurz\\
   University of Bayreuth, 95440 Bayreuth, Germany\\
   sascha.kurz@uni-bayreuth.de}
\date{}
\maketitle   
%%\address{Sascha Kurz, Mathematisches Institut, Universit\"at Bayreuth, 95440 Bayreuth, Germany. E-mail: sascha.kurz@uni-bayreuth.de, 
%%     Phone: +49\,921\,55\,7353, Fax: +49\,921\,55\,7353, Homepage: http://www.wm.uni-bayreuth.de/index.php?id=sk}
\begin{abstract}
\noindent
In 1996 Dan Felsenthal and Mosh\'e Machover considered the following model. An assembly consisting of 
$n$ voters exercises roll-call. All $n!$ possible orders in which the voters may be called are assumed 
to be equiprobable. The votes of each voter are independent with expectation $0<p<1$ for an individual 
vote {\lq\lq}yea{\rq\rq}. For a given decision rule $v$ the \emph{pivotal} voter in a roll-call is the one 
whose vote finally decides the aggregated outcome. It turned out that the probability to be pivotal 
is equivalent to the Shapley-Shubik index. Here we give an easy combinatorial proof of this coincidence, 
further weaken the assumptions of the underlying model, and study generalizations to the case of more than 
two alternatives. 

\medskip


\end{abstract}
%\maketitle                  

\section{Introduction}
\label{sec_introduction}

\noindent
Consider a course in which there are two in-course assessments and an end-of-course 
examination. Assuming the $5$-letter grade system from the US, what should be the course 
result of a student achieving a D, a C, and a B in the three assessments? In practice the most 
common procedure is the following: The letter grades are first converted to numerical values 
and a weighted mean is computed, which then is rounded and converted back to a letter grade. 
There may be additional rules describing the special case of a failure. Two Fs may force a 
final F or an F in the end-of-course examination has to result in an F unless it is balanced 
by two assessments which are graded C or better.

In Germany the school system, roughly spoken, splits into three different branches, where the 
most reputable branch is called {\lq\lq}Gymnasium{\rq\rq}. In Bavaria there are some strict rules 
regulating the access to a Gymnasium. After the forth year the grads of mathematics, german, and 
local history and geography, ranging from 1 to 6, are considered. Iff the sum is at most 7, then 
the student is allowed to access a Gymnasium without any additional acceptance test.

The underlying structure of both examples can be formalized as follows. Let $v\colon J^n\rightarrow K$ 
be a mapping that aggregates $n$ inputs, whose values are contained in $J$, into a single output, which 
is contained in $K$. In our first example we have $n=3$ and $J=K=\{\text{A},\text{B},\text{C},\text{D},\text{F}\}$ 
with $\text{A}>\text{B}>\text{C}>\text{D}>\text{E}$. In our second example we have $n=3$, $J=\{1,2,3,4,5,6\}$, and 
$K=\{\text{denial}, \text{access}\}$ with $\text{access}>\text{denial}$ and $1>2>3>4>5>6$. In the following we 
will assume that the sets $J$ and $K$ are strictly ordered. W.l.o.g.\ we will mostly use the sets $J=\{1,\dots,j\}$ 
and $K=\{1,\dots,k\}$ for two positive integers $j$,  $k$ with the usual ordering over the integers. The school 
context is by far not the only area of application for these aggregation functions. Every committee that takes its 
decision by votes according to some specific voting rule is an example. The binary case, where $J=K=\{\text{yes},\text{no}\}$, 
is extensively treated in the voting literature. For examples with $(j,k)\neq(2,2)$ we refer the interested reader e.g.\ 
to \cite{freixas2009anonymous}.

A classical question in this context asks for the \emph{influence} of a committee member (or voter) on the aggregated 
decision. To this end so-called power indices where introduced. For the binary case, the Shapley-Shubik index, introduced in 
\cite{shapley1954method}, is one of the most commonly used power indices. Besides an axiomatic foundation of 
the Shapley-Shubik index \cite{dubey1975uniqueness}, there is also a picturesque description: Assume that the voters express their support for 
a proposal one after the other. At some point the support is large enough so that the aggregated group decision will 
be an acceptance in any case. The corresponding voter is called a \emph{pivot}. With this, the Shapley-Shubik index of 
a voter $i$ is the fraction of the arrangements of the voters where voter~$i$ is a pivot. In \cite{felsenthal1996alternative} 
the model is extended to a roll-call where each voter can either say {\lq\lq}yes{\rq\rq} or {\lq\lq}no{\lq\lq}. The average 
number of cases where a voter is pivotal coincides with the Shapley-Shubik index. Here the expectation for a {\lq\lq}yes{\rq\rq} 
need not be equal to $\frac{1}{2}$ to preserve this property. A sufficient condition is the independence and equality of 
expectations of the random variables for the voters, as already observed in \cite{mann1964priori} without proof. The proof in
\cite{felsenthal1996alternative} uses the axiomatization of the Shapley-Shubik index as a indirect approach and reports of 
combinatorial difficulties for the direct approach. Here we give an easy combinatorial proof and weaken the assumptions, i.e., 
we assume that the probability for $s$ {\lq\lq}yes{\rq\rq}- and $n-s$ {\lq\lq}no{\rq\rq}-votes only depends on the number $s$. 
Note that the same result was previously obtained in \cite[Proposition 4]{hu2006asymmetric}, as we found out recently.

The generalization of the Shapley-Shubik index to the non-binary case can be traced back at least to \cite{felsenthal1997ternary}, 
where a third input alternative was considered. In full generality this was treated in \cite{freixas2005shapley}, see also 
\cite{tchantcho2008voters}. Here we present a consistent theory that covers the binary case, the general non-binary case, and 
the limiting case with an infinite number of alternatives in both the input and the output. The first basic building blocks
for such a unified theory where sketched in \cite{kurz2014measuring}, see also \cite{kurz2018influence}. 

A main assumption of our considerations is the strict ordering of the sets of alternatives. For the case of unordered alternatives 
we refer the interested reader to e.g.\ \cite{bolger1993value}.

The remaining part of this paper is structures as follows. In Section~\ref{sec_preliminaries} we define the class of games 
with several alternatives in the input and output. The Shapley-Shubik index for simple games is reinterpreted as 
a measure for uncertainty reduction in the roll-call model and generalized to the previously defined more general class of games in 
Section~\ref{sec_reduction_of_uncertainty}. The main result that the chance for being the pivotal player in a simple 
game is almost independent of the probability distribution of the votes is formalized in Theorem~\ref{thm_main} in 
Section~\ref{sec_main_result}. The corresponding proof is purely combinatorial. See \cite[Proposition 4]{hu2006asymmetric} 
for an earlier proof and \cite[Theorem 1]{kurznapelcharacterization} for a characterization of such probability distributions. As a possible 
justification of the presented generalized influence measure we present some first preliminary results on an axiomatization in 
Section~\ref{sec_axioms}. Games with a continuous infinite number of alternatives in the input and output are addressed in 
Section~\ref{sec_limiting_case}. We draw a conclusion in Section~\ref{sec_conclusion}. 

\section{Preliminaries}
\label{sec_preliminaries}

\begin{Definition}
  For $J,K\subseteq \mathbb{R}$ and $n\in\mathbb{N}_{>0}$ the mapping 
  $v\colon J^n\rightarrow  K$ is called \emph{$(J,K)$ game} on $n$ players. For the 
  special sets $J=\{1,\dots,j\}$ and $K=\{1,\dots,k\}$, where $j=|J|$ and $k=|K|$, we 
  speak of \emph{$(j,k)$ games} and denote the set of of $(j,k)$ games on $n$ players by 
  $\mathcal{G}_{j,k}^n$.
\end{Definition}

Since those mappings are not very interesting for $|K|=1$ or $|J|=1$, we assume 
$|J|,|K|,j,k\ge 2$ in the remaining part of the paper. We also speak of \emph{non-trivial} 
games in order to highlight the assumption.  

\begin{Definition}
  A $(J,K)$ game $v$ on $n$ players is called \emph{monotonic} if we have $v(a)\ge v(b)$ for all 
  $a,b\in J^n$ with $a\ge b$, where $(a_1,\dots,a_n)\ge (b_1,\dots,b_n)$ iff $a_i\ge b_i$ for all 
  $1\le i\le n$. By $\mathcal{S}_{j,k}^n$ we denote the set of of all monotonic $(j,k)$ games on $n$ 
  players, which are surjective.\footnote{In \cite{freixas2005shapley} these objects were called 
  $(j,k)$ simple games related to the notion of simple games, where the inputs are labeled in the 
  reverse order. However, in the older papers {\lq\lq}simple{\rq\rq} just refers to the binary input 
  and output. Being more precise, some authors speak of monotonic simple games.} 
\end{Definition}

We remark that $\mathcal{S}_{2,2}^n$ is in bijection with the set of \emph{simple games} on $n$ players, 
see e.g.\ \cite{0943.91005} for an overview on simple games. Simple games are models for simple 
voting situations, where each player can either vote {\lq\lq}yes{\rq\rq} or {\lq\lq}no{\rq\rq}. 
The set $S\subseteq\{1,\dots,n\}:=N$ of {\lq\lq}yes{\rq\rq}-players is called \emph{coalition}. The most 
common formalization of a simple game is given by $\tilde{v}\colon 2^N\rightarrow\{0,1\}$ with 
$\tilde{v}(\emptyset)=0$, $\tilde{v}(N)=1$, and $\tilde{v}(S)\le \tilde{v}(T)$ for all 
$\emptyset\subseteq S\subseteq T\subseteq N$. $\tilde{v}(S)=1$ means an aggregated {\lq\lq}yes{\lq\lq} and 
$\tilde{v}(S)=0$ an aggregated {\lq\lq}no{\rq\rq}. The set notation is just an abbreviation to state that 
the players in $S$ are voting {\lq\lq}yes{\rq\rq} (or $1$) and the players in $N\backslash S$ are voting
{\lq\lq}no{\rq\rq} (or $0$). The bijection from $2^N$ to $\{0,1\}^n$ is just the characteristic vector of 
a set. It is easy to check that a simple game on $n$ players is equivalent to a $(\{0,1\},\{0,1\})$ game.
In the corresponding $(2,2)$ game the inputs and outputs are just increased by $1$. We remark that in the 
field of threshold logic the representation of a (weighted) simple game with $J=K=\{-1,1\}$ is more common.

\begin{Definition}
  Two players $1\le i,h\le n$ are called \emph{equivalent} in a $(J,K)$ game $v$ on $n$ players, if we have 
  $v(a)=v(\pi(a))$ for all $a\in J^n$, where $\pi$ is the transposition between $i$ and $h$. 
\end{Definition}

A classical question in this context asks for the \emph{influence} of a committee member (or voter) on the aggregated 
decision. For a simple game $\tilde{v}\colon 2^N\rightarrow\{0,1\}$ the so-called \emph{Shapley-Shubik index}, see
\cite{shapley1954method}, of player $1\le i\le n$ in $\tilde{v}$ is given by
\begin{equation}
  \label{eq_ssi}
  \tilde{\varphi}_i(\tilde{v})=\frac{1}{n!}\cdot \sum_{\emptyset \subseteq S\subseteq N\backslash\{i\}} 
  |S|!\cdot (n-1-|S|)!\cdot\left(v(S\cup \{i\})-v(S)\right).
\end{equation}
It can be interpreted as a weighted marginal contribution and is just a specialization of the 
Shapley value for transferable utility (TU) games, see \cite{shapley_value}. In the next section we will 
describe an influence measure for $(j,k)$ games similar to the Shapley-Shubik index.

The extreme case of having absolutely no influence is captured by:
\begin{Definition}
  \label{def_null_player}
  Let $v$ be a $(J,K)$ game on $n$ players. If we have 
  $$v(a_1,\dots,a_n)=v(a_1,\dots,a_{i-1},x,a_{i+1},\dots,a_n)$$ 
  for all $a_1,\dots,a_n,x\in J$, 
  player $i$ is called a \emph{null player}.
\end{Definition}
Of course this general definition is consistent with the definition of a null player in a simple game 
$\tilde{v}$ and we have $\tilde{\varphi}_i(\tilde{v})=0$ for each null player~$i$. For equivalent players 
$h,i$ we have $\tilde{\varphi}_h(\tilde{v})=\tilde{\varphi}_i(\tilde{v})$.

\section{The Shapley-Shubik index as a measurement for the reduction of uncertainty} 
\label{sec_reduction_of_uncertainty}

\noindent
If the votes of all $n$ players are known, then the aggregated decision is uniquely determined by the 
underlying game. In order to evaluate the \emph{influence} of each player on the final outcome, one can 
think of the voting situation as a roll-call, c.f.~\cite{felsenthal1996alternative}. Each player declares
her vote one after the other. In the case of a binary outcome, i.e., $|K|=2$, there exists a certain 
player $i$ whose declaration uniquely determines the outcome for the first time. The respective player 
is commonly called \emph{pivotal}, cf.~Definition~\ref{def_pivot}. Of course the pivotality depends on the ordering of the players in 
this context. For more than two output alternatives the set of possible outcomes may shrink several times.

\begin{Definition}
  Let $v$ be a $(J,K)$ game on $n$ players and $S_n$ be the set of all bijections of $N$, i.e., the set 
  of all permutations of $N=\{1,\dots,n\}$. For each $1\le h\le n$ we set 
  $\tau_h\colon \mathcal{G}^n_{J,K}\times S_n\times J^n\rightarrow \{0,\dots,|K|-1\}$, 
  \begin{eqnarray*}
    \tau_h(v,\pi,a_1,\dots,a_n)&=&
    \left|\left\{ v(a_1',\dots,a_n')\,:\, a_l'=a_l\text{ if } \pi(l)< h \text{ and } a_l'\in J\text{ otherwise}\right\}\right|\\
    &&-\left|\left\{ v(a_1',\dots,,a_n')\,:\, a_l'=a_l\text{ if } \pi(l)\le h \text{ and } a_l'\in J\text{ otherwise}\right\}\right|,
  \end{eqnarray*}
  where $\mathcal{G}^n_{J,K}$ denotes the set of $(J,K)$ game on $n$ players. 
\end{Definition}

So, given an ordering $\pi$ of the $n$ players and a specific input $(a_1,\dots,a_n)$, the 
value of $\tau_h(v,\pi,a_1,\dots,a_n)$ displays the decrease of our uncertainty of the final outcome after 
the $h$th player, according to $\pi$, has declared her vote. The initial uncertainty is $|\operatorname{im}(v)|-1$, 
where 
$$
  \operatorname{im}(v)=\highlight{\left\{v(a_1,\dots,a_n)\,:\,a_i\in J\, \forall 1\le i\le n\right\}}
$$
denotes the image of the mapping $v$, i.e., there are $|\operatorname{im}(v)|$ possible outcomes before the first 
player declares her vote and a unique outcome after the last player has declared her vote. By summing over all possible 
orderings and all possible inputs we obtain a general measurement for \emph{influence} after normalization:
\begin{Definition}
  \label{def_influence_measure}
  For each positive integer $n$, each integer $1\le i\le n$, and each sets $J,K\subseteq\mathbb{R}$ with $|J|,|K|\ge 2$ we set 
  $\varphi_i\colon \mathcal{G}_{J,K}^n\rightarrow\mathbb{R}_{\ge 0}$ with 
  \begin{equation}
    \label{eq_influence_measure}
    \varphi_i(v)=\frac{1}{n!}\cdot\frac{1}{|J|^n}\cdot \frac{1}{|\operatorname{im}(v)|-1}\cdot
    \sum_{(\pi,a_1,\dots,a_n)\in S_n\times J^n} \tau_{\pi^{-1}(i)} (v,\pi,a_1,\dots,a_n)
  \end{equation}
  if $|\operatorname{im}(v)|>1$ and $\varphi_i(v)=0$ otherwise. 
\end{Definition} 

We remark that we have $|\operatorname{im}(v)|=|K|$ if $v$ is surjective. Via the normalization factors, the \emph{influence 
measure} of Definition~\ref{def_influence_measure} obtains a nice property that is called \emph{efficiency} in the context 
of power indices.

\begin{Lemma}
  \label{lemma_efficient}
  For each $v\in\mathcal{G}_{J,K}^n$ with $|\operatorname{im}(v)|>1$ we have 
  $\sum\limits_{i=1}^n \varphi_i(v)=1$ and $\varphi_i(v)\in[0,1]$ for all $1\le i\le n$.\footnote{We remark 
  that the Shapley-Shubik index for simple games, based on Equation~(\ref{eq_ssi}), is not efficient in all 
  non-monotonic cases. For $\tilde{v}(\emptyset)\!=\!v(\{2\})\!=\!v(\{1,2\})\!=\!0$, $\tilde{v}(\{1\})\!=\!1$ we would have 
  $\tilde{\varphi}_1(\tilde{v})=\highlight{\frac{1}{2}}$ and $\tilde{\varphi}_2(\tilde{v})=\highlight{-\frac{1}{2}}$.} 
\end{Lemma} 
\begin{Proof}
  For each $\pi\in S_n$ and each $(a_1,\dots,a_n)\in J^n$ we have 
  $$
    \sum_{i=1}^n \tau_{\pi^{-1}(i)}(v,\pi,a_1,\dots,a_n)=
    \sum_{h=1}^n \tau_h(v,\pi,a_1,\dots,a_n)= |\operatorname{im}(v)|-1.
  $$
\end{Proof}
  
\begin{Definition}
  \label{def_output_rough}
  We call $v\in\mathcal{G}_{j,k}^n$ \emph{output-rough} if
  $$
    \left\{ v(a_1',\dots,a_n')\,:\, a_l'=a_l\text{ if } \pi(l)\le h \text{ and } a_l'\in J\text{ otherwise}\right\}
  $$   
  is an interval\footnote{We call a set $M\subseteq \mathbb{Z}$ an interval if there exist $a,b\in \mathbb{R}$ with 
  $M=[a,b]\cap\mathbb{Z}$.} for all $(a_1,\dots,a_n)\in J^n$, $\pi \in S_n$, $1\le h\le n$.
\end{Definition}

We remark each $(j,2)$ game is output-rough. An example of a
surjective and monotonic $(2,3)$ game that is not output-rough ins given in Example~\ref{ex_2}. An example 
of a surjective, monotonic, output-rough $(2,3)$ game $v$ is given by $v(1,1)=1$, $v(2,2)=3$, and $v(1,2)=v(2,1)=2$.

For the case of surjective, monotonic, output-rough $(j,k)$ games with $j,k\ge 2$, Equation~(\ref{eq_influence_measure}) 
can be simplified significantly.

\begin{Lemma}
  \label{lemma_influence_measure_monotonic}
  For each non-trivial, surjective, monotonic, output-rough $(j,k)$ game $v$ on $n$ players and each integer $1\le i\le n$ we have
  \begin{eqnarray}
    \label{eq_influence_measure_monotonic}
    \varphi_i(v)&=&\frac{1}{n!}\cdot\frac{1}{j^n}\cdot \frac{1}{k-1}\cdot\\
    &&\!\!\!\!\!\!\!\!\!\!\!\!\!\!\!\!\!\!\!\!\sum_{(\pi,a=(a_1,\dots,a_n))\in S_n\times J^n} \!\!\!\!\!\!\!\!\!\!\!\!\!\!\!\!\!\!\!\!\!\! 
     \left(\overline{v}_{\pi,\pi^{-1}(i)\!-\!1}(a)-\underline{v}_{\pi,\pi^{-1}(i)\!-\!1}(a)\right)
     -\left(\overline{v}_{\pi,\pi^{-1}(i)}(a)-\underline{v}_{\pi,\pi^{-1}(i)}(a)\right),\nonumber
  \end{eqnarray}  
  where $\overline{v}_{\pi,h}(a)=v(a_1',\dots,a_n')$ with $a_l'=a_l$ if $\pi(l)\le h$ and $a_l'=j$ otherwise;
  $\underline{v}_h(b)=v(b_1',\dots,b_n')$ with $b_l'=b_l$ if $\pi(l)\le h$ and $b_l'=1$ otherwise. 
\end{Lemma} 
\begin{Proof}
  Due to monotonicity and output-roughness we have 
  $$
    \left\{ v(a_1',\dots,a_n')\,:\, a_l'=a_l\text{ if }\pi(l)\le h \text{ and } a_l'\in J\text{ otherwise}\right\}=
    \left\{\underline{v}_{\pi,h}(a),\dots,\overline{v}_{\pi,h}(a)\right\}
  $$
  for all $1\le h\le n$ and all $a=(a_1,\dots,a_n)\in J^n$.
\end{Proof}

For $j=k=2$ the formula for $\varphi_i$ can be simplified significantly, which is the topic of the next section. Some 
more notation can be introduced for the slightly more general case $k=2$.

\begin{Definition}
  \label{def_pivot}
  For $j\ge 2$, let $v$ be a surjective $(j,2)$ game on $n$ players. Given an input vector $a\in J	^n$ and 
  an ordering $\pi\in S_n$ we call a player $1\le i\le n$ a \emph{pivot} for $a$, $\pi$ in $v$ if 
  $\tau_{\pi^{-1}(i)}(v,\pi,a)=1$.
\end{Definition}

In other words, a player $i$ is a pivot, if the declarations before player $i$ still allow both elements of $K$ as 
possible outcomes while player $i$ fixes the outcome (to either $1$ or $2$). For a non-trivial, surjective, monotonic 
$(j,k)$ game $v$ we may describe 
$$
  M=\left\{ v(a_1',\dots,a_n')\,:\, a_l'=a_l\text{ if } \pi(l)\le  \pi(i) \text{ and } a_l'\in J\text{ otherwise}\right\}
$$
by $b_1^M,\dots,b_{k-1}^M\in\{-1,0,1\}$ with 
\begin{itemize}
  \item $b_h^M=-1$ iff $m\le h$ for all $m\in M$;
  \item $b_h^M=1$ iff $m> h$ for all $m\in M$;
  \item $b_h^M=0$ otherwise.
\end{itemize} 
If $b_h^M$ switches from a zero to a non-zero value, then we may call player $i$ an $h$-pivot for $a$, $\pi$ in $v$. 
The number of pivots for given $a$, $\pi$ is $k-1$. With this, $\varphi_i(v)$ is equal to the probability of a player $i$ 
to be pivotal, assuming equiprobable input vectors and orderings.

The notion of minimal winning coalitions can be generalized to non-trivial, surjective, monotonic $(j,2)$ games.

\begin{Definition}
  Let $v$ be a non-trivial, surjective, monotonic $(j,2)$ game on $n$ players. A vector $a\in J^n$ is called a
  \emph{winning vector} if $v(a)=2$ and \emph{losing vector} otherwise. If $a$ is a winning vector, but all vectors 
  $a'<a$ are losing, then $a$ is called a \emph{minimal winning vector}. Analogously a losing vector $a$ is called 
  \emph{maximal losing vector} if $a'$ is a winning vector for all $a'>a$.\footnote{We write $a>b$ for $a\ge b$ and $a\neq b$.}
\end{Definition}
We remark that $v$ is uniquely described by its set of minimal winning or its set of maximal losing vectors. The minimal 
winning vectors of Example~\ref{ex_1} are given by $(2,3)$ and $(3,2)$. The corresponding maximal losing vectors are given 
by $(1,3)$, $(3,1)$, and $(2,2)$.

\section{The main result}
\label{sec_main_result}

\noindent
Now we focus on the special case $j=k=2$ and relate Equation~(\ref{eq_ssi}) with 
Equation~(\ref{eq_influence_measure_monotonic}). To this end we have to introduce 
some further notions for simple games.

\begin{Definition} 
  Let $\tilde{v}$ be a simple game on $n$ players. A coalition $S\subseteq N\backslash\{i\}$ 
  is called an \emph{$i$-swing} if $v(S\cup\{i\})-v(S)=1$, i.e., $v(S\cup\{i\})=1$ and $v(S)=0$. 
\end{Definition}

With this we can rewrite Equation~(\ref{eq_ssi}) to
$$
  \tilde{\varphi}_i(\tilde{v})=\frac{1}{n!}\cdot \sum_{S\text{ is an }i-\text{swing}} 
  |S|!\cdot (n-1-|S|)!.
$$
Let $v\in\mathcal{S}_{2,2}^n$ correspond to $\tilde{v}$ and $S$ be an $i$-swing 
for an arbitrary but fixed player~$i$. We set $a_l=2$ if $l\in S$ and $a_l=1$ otherwise for all $l\in N\backslash\{i\}$. 
With this we have $v(a_1,\dots,a_{i-1},1,a_{i+1},\dots,a_n)=1$ and $v(a_1,\dots,a_{i-1},2,a_{i+1},\dots,a_n)=2$. So, 
if the players in $S$ are asked first, then still both outcomes $1$ and $2$ are possible. If player~$i$ says $2$, then 
the outcome is fixed to $2$ due to the monotonicity of $v$. There are exactly $|S|!\cdot (n-1-|S|)!$ orderings 
where the set of players before player~$i$ coincides with $S$. There is another interpretation: If the players of 
$N\backslash (S\cup\{i\})$ are asked first, then still both outcomes $1$ and $2$ are possible. If player~$i$ says $1$, then 
the outcome is fixed to $1$ due to the monotonicity of $v$. From these observations we conclude:

\begin{Lemma}
  \label{lemma_ssi_extreme_cases}
  Let $\tilde{v}$ be a simple game on $n$ players and $v$ be the corresponding $(2,2)$ game.
  We have 
  \begin{eqnarray*}
  \tilde{\varphi}_i(\tilde{v}) &=&
    \frac{1}{n!}\cdot \sum_{S\text{ is an }i-\text{swing}} |S|!\cdot (n-1-|S|)!\\
    &=&
    \frac{1}{n!}\cdot
    \sum_{(\pi,a=(2,\dots,2))\in S_n\times J^n} \!\!\!\!\!\!\!\!\!\!\!\! 
     \left(\overline{v}_{\pi,\pi^{-1}(i)\!-\!1}(a)\!\!-\!\!\underline{v}_{\pi,\pi^{-1}(i)\!-\!1}(a)\right)
     -\left(\overline{v}_{\pi,\pi^{-1}(i)}(a)\!\!-\!\!\underline{v}_{\pi,\pi^{-1}(i)}(a)\right)\\
    &=&
    \frac{1}{n!}\cdot
    \sum_{(\pi,a=(1,\dots,1))\in S_n\times J^n} \!\!\!\!\!\!\!\!\!\!\!\! 
     \left(\overline{v}_{\pi,\pi^{-1}(i)\!-\!1}(a)\!\!-\!\!\underline{v}_{\pi,\pi^{-1}(i)\!-\!1}(a)\right)
     -\left(\overline{v}_{\pi,\pi^{-1}(i)}(a)\!\!-\!\!\underline{v}_{\pi,\pi^{-1}(i)}(a)\right).
  \end{eqnarray*}  
\end{Lemma}  

So, the Shapley-Shubik index is the same as the influence measure from Definition~\ref{def_influence_measure} 
when all players say $2$ ({\lq\lq}yes{\rq\rq}) or all players say $1$ ({\lq\lq}no{\rq\rq}) in the roll-call 
model. Equation~(\ref{eq_ssi}) can be seen as a computational simplification. We will see shortly that the connection 
between Equation~(\ref{eq_ssi}) and the roll-call model is far more general than suggested by Lemma~\ref{lemma_ssi_extreme_cases}.  
To this end denote by $\chi_h(a)$ the number of $a_i$ which are equal to $h$.

\begin{Lemma}
  \label{lemma_contribution}
  For $v\in\mathcal{S}_{2,2}^n$ we have 
  \begin{equation}
    \label{eq_contribution}
    \left(\overline{v}_{\pi,\pi^{-1}(i)\!-\!1}(a)-\underline{v}_{\pi,\pi^{-1}(i)\!-\!1}(a)\right)
     -\left(\overline{v}_{\pi,\pi^{-1}(i)}(a)-\underline{v}_{\pi,\pi^{-1}(i)}(a)\right)
     =1
  \end{equation}
  iff either
  $$
    \underline{v}_{\pi,\pi^{-1}(i)\!-\!1}(a)=\underline{v}_{\pi,\pi^{-1}(i)}(a)= 
    \overline{v}_{\pi,\pi^{-1}(i)}(a)=1,\, \overline{v}_{\pi,\pi^{-1}(i)\!-\!1}(a)=2,\, a_i=1    
  $$
  or
  $$
    \underline{v}_{\pi,\pi^{-1}(i)\!-\!1}(a)=1,\, \underline{v}_{\pi,\pi^{-1}(i)}(a)= 
    \overline{v}_{\pi,\pi^{-1}(i)}(a)=\overline{v}_{\pi,\pi^{-1}(i)\!-\!1}(a)=2,\,a_i=2.
  $$
\end{Lemma}
\begin{Proof}
  Due to monotonicity the following cases are possible:
  \begin{center}
    \begin{tabular}{rrrr}
      \hline
      $\underline{v}_{\pi,\pi^{-1}(i)\!-\!1}(a)$ & $\underline{v}_{\pi,\pi^{-1}(i)}(a)$ & 
      $\overline{v}_{\pi,\pi^{-1}(i)}(a)$ &  $\overline{v}_{\pi,\pi^{-1}(i)\!-\!1}(a)$\\
      \hline 
      1 & 1 & 1 & 1 \\
      1 & 1 & 1 & 2 \\
      1 & 1 & 2 & 2 \\
      1 & 2 & 2 & 2 \\  
      2 & 2 & 2 & 2 \\
      \hline
    \end{tabular}
  \end{center}
\end{Proof}

\begin{Lemma}
  \label{lemma_binom_sum_identity}
  For integers $0\le s\le h-1\le n-1$ we have
  $$
    \sum_{l=0}^{n-h} {{s+l}\choose s}\cdot {{n-s-1-l}\choose{h-s-1}}={n\choose h}.
  $$
\end{Lemma}
\begin{Proof}
  The stated summation formula can be concluded from Vandermonde's Identity.  
  %http://www.trans4mind.com/personal_development/mathematics/series/summingBinomialCoefficients.htm#Identity_5:_Minus_k_in_one_upper_index
  Here we give a direct combinatorial proof by double counting. The number of ways to choose $h$ out of $n$ objects is 
  given by ${n \choose h}$. For a selection let $l$ be an integer such that the $(s+1)$th chosen object is labeled 
  $s+l+1$, where we assume labels from $1$ to $n$. Here $l$ can range from $0$ to $n-h$ and is uniquely determined. Since exactly 
  $s$ elements have to be chosen before the $(s+1)$th element and $h-s-1$ elements have to be chosen after the $(s+1)$th element, there
  are ${{s+l}\choose s}\cdot{{n-s-1-l}\choose{h-s-1}}$ possibilities. Summing over the possible values for $l$ gives the stated formula.      
\end{Proof}
  
\begin{Lemma}
  \label{lemma_combinatorial_key_result}
  Let $\tilde{v}$ be a simple game on $n$ players and $v$ be the corresponding $(2,2)$ game. For each 
  $1\le i\le n$ and each $0\le h\le n$ we have
  \begin{eqnarray*}
    &&{n \choose h}\cdot \sum_{S\text{ is an }i-\text{swing}} |S|!\cdot (n-1-|S|)!\\
    &=&
    \sum_{(\pi,a)\in S_n\times J^n\text{ with }\chi_2(a)=h} \!\!\!\!\!\!\!\!\!\!\!\! 
     \left(\overline{v}_{\pi,\pi^{-1}(i)\!-\!1}(a)-\underline{v}_{\pi,\pi^{-1}(i)\!-\!1}(a)\right)
     -\left(\overline{v}_{\pi,\pi^{-1}(i)}(a)-\underline{v}_{\pi,\pi^{-1}(i)}(a)\right).
  \end{eqnarray*}   
\end{Lemma} 
\begin{Proof}
  We prove the stated equation by double counting. The right hand side clearly counts the number of cases 
  $(\pi,a)\in S_n\times J^n$ with $\chi_2(a)=h$ where Equation~(\ref{eq_contribution}) is satisfied.
  
  Now let $\pi\in S_n$ and $a\in J^n$ with $\chi_2(a)=h$, where Equation~(\ref{eq_contribution}) is satisfied, 
  be arbitrary but fixed. Set $F=\{l\in N\,:\, \pi(l)<\pi(i)\}$ and 
  $B=\{l\in N\,:\, \pi(l)>\pi(i)\}$, so that $F\cup\{i\}\cup B=N$ is a partition. If 
  $a_i=2$, then we set $S=\{l\in F\,:\, a_l=2\}$. If $a_i=1$, then we choose 
  $S\subseteq N\backslash\{i\}$ such that $\{l\in F\,:\,a_l=1\}=N\backslash (S\cup\{i\})$. 
  We can easily check that $S$ is an $i$-swing in both cases. 
  
  Now we start with a fixed $i$-swing $S$ and count the corresponding pairs $(\pi,a)$ as described above. We can 
  easily check that exactly one of the cases $|S|\le h-1$ or $|N\backslash(S\cup\{i\})|\le n-h-1$ is satisfied.
  \begin{itemize}
    \item For the second case of Lemma~\ref{lemma_contribution} we need $|S|\le h-1$ and $a_m=2$ for all $m\in S\cup\{i\}$.
          Exactly $h-s-1$ out of the remaining $n-s-1$ players have to vote $2$ so that $\chi_2(a)=h$. Choose a integer $l$ 
          so that the set $F$, as described above, has cardinality $|S|+l$. Thus, we have $|B|=n-|S|-1-l$, $0\le l\le n-h$ 
          and $l$ out of the $n-h$ players voting $1$ have to be chosen for $F$. Since there are $|F|!\cdot |B|!$ fitting
          permutations $\pi$, we obtain
          \begin{eqnarray*}
            && {{n-s-1}\choose{h-s-1}}\cdot\sum_{l=0}^{n-h} (s+l)!\cdot (n-s-1-l)!\cdot {{n-h}\choose l}\\
            &=& \sum_{l=0}^{n-\highlight{h}} \frac{(n-s-1)!\cdot(s+l)!\cdot(n-s-1-l)!\cdot(n-h)!}{(h-s-1)!\cdot(n-h)!\cdot l!\cdot(n-h-l)!}\\
            &=& s!\cdot(n-s-1)!\cdot \sum_{l=0}^{n-h} {{s+l}\choose s}\cdot{{n-s-1-l}\choose {h-s-1}}\\
            &\overset{\text{Lemma~\ref{lemma_binom_sum_identity}}}{=}& {n\choose h}\cdot s!\cdot(n-s-1)!
          \end{eqnarray*}   
          cases, where we use $|S|=s$ as abbreviation.
    \item For the first case of Lemma~\ref{lemma_contribution} we need $|N\backslash(S\cup\{i\})|\le n-h-1$ and 
          $a_m=1$ for all $m\in N\backslash S$. Exactly $s-h$ out of the remaining $s$ players have to vote $1$ so 
          that $\chi_2(a)=h$. Choose a integer $l$ so that the set $F$, as described above, has cardinality $n-s-1+l$. 
          Thus, we have $|B|=s-l$, $0\le l\le h$ and $l$ out of the $h$ players voting $2$ have to be chosen for $F$. 
          Since there are $|F|!\cdot |B|!$ fitting permutations $\pi$, we obtain
          \begin{eqnarray*}
            && {{s}\choose{s-h}}\cdot\sum_{l=0}^{h} (n-s-1+l)!\cdot (s-l)!\cdot {{h}\choose l}\\
            &=& s!\cdot(n-s-1)!\cdot \sum_{l=0}^{h} {{s'+l}\choose {s'}}\cdot {{n-s'-1-l}\choose {h'-s'-1}}\\
            &=& \highlight{ s!\cdot(n-s-1)!\cdot \sum_{l=0}^{n-h'} {{n-s-1+l}\choose {n-s-1}}\cdot {{s-l}\choose {s-h}} }
            %%&\overset{\text{Lemma~\ref{lemma_binom_sum_identity}}}{=}& {n\choose n-h}\cdot s!\cdot(n-s-1)!
            %%={n\choose h}\cdot s!\cdot(n-s-1)!
          \end{eqnarray*}          
          cases \highlight{using $h'=n-h$ and $s'=n-s-1$. Since $n-s-1\le n-h-1$ we have $0\le s'\le h'-1\le n-1$, so that we can 
          apply Lemma~\ref{lemma_binom_sum_identity}. Using ${n\choose {h'}}={n\choose {n-h}}={n\choose {h}}$, we have
          ${n\choose h}\cdot s!\cdot(n-s-1)!$ cases.} 
  \end{itemize}  
\end{Proof}

Applying Lemma~\ref{lemma_combinatorial_key_result} yields our main result:

\begin{Theorem} (Cf.~\cite[Proposition 4]{hu2006asymmetric})
  \label{thm_main}
  Let $\tilde{v}$ be a simple game on $n$ players, $v$ be the corresponding $(2,2)$ game, and $p\colon J^n\rightarrow [0,1]$ 
  a probability measure with $p(a)=p(b)$ for all $a,b\in J^n$ with $\chi_2(a)=\chi_2(b)$. For each $1\le i\le n$ we have
  \begin{equation*}
    \label{eq_main}
    \tilde{\varphi}_i(\tilde{v})=
    \frac{1}{n!}\cdot\!\!\!\!\!\!\!\!\!\!\!\!\!\!\!\!\!\!\!\!
    \sum_{(\pi,a=(a_1,\dots,a_n))\in S_n\times J^n} \!\!\!\!\!\!\!\!\!\!\!\!\!\!\!\!\!\!\! p(a)\cdot\Big(  
     \left(\overline{v}_{\pi,\pi^{-1}(i)\!-\!1}(a)-\underline{v}_{\pi,\pi^{-1}(i)\!-\!1}(a)\right)
     -\left(\overline{v}_{\pi,\pi^{-1}(i)}(a)-\underline{v}_{\pi,\pi^{-1}(i)}(a)\right)\Big). 
  \end{equation*} 
\end{Theorem}

We remark that for $p(a)=q^{\chi_2(a)}\cdot(1-q)^{n-\chi_2(a)}$ we obtain the result that the Shapley-Shubik 
index of a player~$i$ in a simple game is equal to the probability of $i$ being pivotal, where the players' votes 
are independent and the individual {\lq\lq}yes{\rq\rq}-votes have an expectation of $0\le q\le 1$, cf.~\cite{mann1964priori}. 

As a refinement we weaken the assumption that the individuals votes are independent to, lets say, \textit{anonymous} 
probabilities\footnote{Cf.~\cite{freixas2012probabilistic}.} 
for the votes, i.e., the probability for a vector of votes only depends on the number of {\lq\lq}yes{\rq\rq}-votes, which seems 
to be a very reasonable assumption. 

For the reverse statement of Theorem~\ref{thm_main} we refer to \cite[Theorem 1]{kurznapelcharacterization}. 

Setting $p(2,\dots,2)=1$ and $p(a)=0$ otherwise, we can obtain Equation~(\ref{eq_ssi}) as a computational simplification 
of the general roll-call model with \textit{anonymous} probabilities for the votes in Theorem~\ref{thm_main}. (If $a_i=2$ for 
all $1\le i\le n$ we can simplify 
$$
  \left(\overline{v}_{\pi,\pi^{-1}(i)\!-\!1}(a)-\underline{v}_{\pi,\pi^{-1}(i)\!-\!1}(a)\right)
  -\left(\overline{v}_{\pi,\pi^{-1}(i)}(a)-\underline{v}_{\pi,\pi^{-1}(i)}(a)\right)
$$ 
to $\underline{v}_{\pi,\pi^{-1}(i)}(a)-\underline{v}_{\pi,\pi^{-1}(i)\!-\!1}(a)$, which is equivalent to $v(S\cup\{i\})-v(S)$ 
for $S=\{l\in N\,:\, \pi(l)<\pi(i)\}$.)

For $p(a)=\frac{1}{2^n}$ we obtain $\tilde{\varphi}_i(\tilde{v})=\varphi_i(v)$, which gives some justification for 
calling the influence measure from Definition~\ref{def_influence_measure} the Shapley-Shubik index for $(j,k)$ games. Using the 
correspondence between $i$-swings and input vectors satisfying Equation~(\ref{eq_contribution}) from the proof 
of Theorem~\ref{thm_main} we can directly prove this equation. Starting from an $i$-swing $S$ we write the set $F$ as $S\cup X$ 
and set $|X|=x$. In contrast to the proof of Theorem~\ref{thm_main} we make no assumption on the values of $a_l$ for $l\in B$. 
So, we have ${{n-s-1}\choose{x}}$ possibilities for $X$, $2^{n-s-x-1}$ possibilities for $a$, and $(s+x)!\cdot (n-s-x-1)!$ possibilities 
for $\pi$, where $0\le x\le n-s-1$. We compute
\begin{eqnarray*}
  &&\sum_{x=0}^{n-s-1} (s+x)!\cdot (n-s-x-1)!\cdot {{n-s-1}\choose{x}} \cdot 2^{n-s-x-1}\\
  &=& s!\cdot (n-s-1)!\cdot \sum_{x=0}^{n-s-1} {{s+x}\choose x} \cdot 2^{n-s-x-1} \\
  &=& s!\cdot (n-s-1)!\cdot\sum_{y=0}^{n-s-1} {{n-1-y}\choose {n-s-1-y}} \cdot 2^{y}
\end{eqnarray*}      
corresponding possibilities.

For the first case of Lemma~\ref{lemma_contribution} we similarly have 
${{s}\choose{x}}$ possibilities for $X$, $2^{s-x}$ possibilities for $a$, and $(n-s-1+x)!\cdot (s-x)!$ possibilities 
for $\pi$, where $0\le x\le s$. We compute
\begin{eqnarray*}
  &&\sum_{x=0}^{s} (s-x)!\cdot (n-s+x-1)!\cdot {{s}\choose{x}} \cdot 2^{s-x}\\
  &=& s!\cdot (n-s-1)!\cdot \sum_{x=0}^{s} {{n-s-1+x}\choose x} \cdot 2^{s-x} \\
  &=& s!\cdot (n-s-1)!\cdot \sum_{y=0}^{s} {{n-1-y}\choose {s-y}} \cdot 2^{y}
\end{eqnarray*}      
corresponding possibilities.

The equation $\tilde{\varphi}_i(\tilde{v})=\varphi_i(v)$, which was indirectly proven in \cite{felsenthal1996alternative} now 
follows from
\begin{equation}
  \label{eq_comb_identity}
  \left(\sum_{y=0}^{n-s-1} {{n-1-y}\choose {n-s-1-y}} \cdot 2^{y}\right)\,+\,
  \left(\sum_{y=0}^{s} {{n-1-y}\choose {s-y}} \cdot 2^{y}\right)=2^n.
\end{equation}
This is indeed a very interesting identity on its own. According to the computer algebra system \texttt{Maple} the first sum 
can be simplified to $2^n-\frac{1}{2}\cdot {n\choose s}\cdot \operatorname{hypergeom}\!\left([1,n+1],,\frac{1}{2}\right)$ and 
the second sum can be simplified to $2^n-\frac{1}{2}\cdot {n\choose {s+1}}\cdot \operatorname{hypergeom}\!\left([1,n+1],[s+2],\frac{1}{2}\right)$.
%% check if both expressions are typed correctly 
In the theory of hypergeometric series, see e.g.~\cite{gasper2004basic}, the corresponding identity
$$
  {n\choose s}\cdot \operatorname{hypergeom}\!\left([1,n+1],,\frac{1}{2}\right)
  +
  {n\choose {s+1}}\cdot \operatorname{hypergeom}\!\left([1,n+1],[s+2],\frac{1}{2}\right)
  =2^{n+1}
$$
%% is Maple able to simplify this expression???
might be well known. Arguably, one can speak of {\lq\lq}rather formidable combinatorial difficulties{\rq\rq}, as done 
in \cite{felsenthal1996alternative}. Plugging in small values of $s$ into Equation~(\ref{eq_comb_identity}) and explicitly 
evaluating the second sum gives some nice explicit identities. For $s=0$ we obtain the well known 
geometric series $\sum_{y=0}^{n-1} 2^{y}=2^n-1$. For $s=1,2$ we obtain
$$
  \sum_{y=0}^{n-2} (n-1-y) \cdot 2^{y}=2^n -n-1  
$$
and
$$
  \sum_{y=0}^{n-3} (n-1-y)\cdot(n-2-y) \cdot 2^{y-1}=2^n -\frac{n^2+n+2}{2}.  
$$

We finish this section by two examples showing that the situation of Theorem~\ref{thm_main} can not be generalized to arbitrary 
parameters $j$ and $k$. 

\begin{Example}
  \label{ex_1}
  Let $v$ be a $(3,2)$ game on $2$ players with $v(a_1,a_2)=2$ iff $(a_1,a_2)\ge(3,2)$ or $(a_1,a_2)\ge (2,3)$. So, the game 
  $v$ is surjective, monotonic, and output-rough. In the following 
  table we list the pivotal player for all combinations of the input vector $a$ and the ordering $\pi$. We assume that each player 
  votes $l$ with probability $p_l$, i.e., $p_1,p_2,p_3\ge 0$ and $p_1+p_2+p_3=1$, and that the votes of the two players are 
  independent.
  \begin{center}
    \begin{tabular}{|rcc|r|}
      \hline
      $a\backslash \pi$ & $(1,2)$ & $(2,1)$ & $p(a)$\\
      \hline
      $(1,1)$ & 1 & 1 & $p_1\cdot p_1$\\
      $(1,2)$ & 1 & 1 & $p_1\cdot p_2$\\
      $(1,3)$ & 1 & 1 & $p_1\cdot p_3$\\
      $(2,1)$ & 2 & 1 & $p_2\cdot p_1$\\
      $(2,2)$ & 2 & 1 & $p_2\cdot p_2$\\
      $(2,3)$ & 2 & 1 & $p_2\cdot p_3$\\
      $(3,1)$ & 1 & 1 & $p_3\cdot p_1$\\
      $(3,2)$ & 1 & 1 & $p_3\cdot p_2$\\
      $(3,3)$ & 1 & 1 & $p_3\cdot p_3$\\
      \hline
    \end{tabular}  
  \end{center}
  Given the probabilities the (generalized) influence measure for player $2$ is given by $\frac{p_1}{2}$ and the (generalized) 
  influence measure for player $1$ is given by $1-\frac{p_1}{2}$, i.e., the values are not independent from the probability 
  distribution. 
\end{Example}

\begin{Example}
  \label{ex_2}
  Let $v$ be a $(2,3)$ game on $2$ players with $v(2,2)=v(2,1)=3$, $v(1,2)=2$, and $v(1,1)=1$. So, the game 
  $v$ is surjective and monotonic but not output-rough. In the following table we list for each player the reduction 
  of uncertainty $\tau$ for all combinations of the input vector $a$ and the ordering $\pi$. We assume that each player  
  votes $l$ with probability $p_l$, i.e., $p_1,p_2\ge 0$ and $p_1+p_2=1$, and that the votes of the two players are 
  independent.
  \begin{center}
    \begin{tabular}{|rcc|r|}
      \hline
      $a\backslash \pi$ & $(1,2)$ & $(2,1)$ & $p(a)$\\
      \hline
      $(1,1)$ & 1:1,2:1 & 1:1,2:1 & $p_1\cdot p_1$\\
      $(1,2)$ & 1:1,2:1 & 1:1,2:1 & $p_1\cdot p_2$\\
      $(2,1)$ & 1:2,2:0 & 1:1,2:1 & $p_2\cdot p_1$\\
      $(2,2)$ & 1:2,2:0 & 1:1,2:1 & $p_2\cdot p_2$\\
      \hline
    \end{tabular}  
  \end{center}
  Given the probabilities the (generalized) influence measure for player $2$ is given by $\frac{1+p_1}{4}$ and the (generalized) 
  influence measure for player $1$ is given by $\frac{3-p_1}{2}$, i.e., the values are not independent from the probability 
  distribution. 
\end{Example}

Also for $(2,2)$ games the conditions of Theorem~\ref{thm_main} can not be weakened too much.

\begin{Example}
  Let $v$ be the $(2,2)$ game on $3$~players with $v(a_1,a_2,a_3)=2$ iff $a_1=1$ and $a_2\neq a_3$. So, $v$ is surjective and 
  output-rough but not monotonic. In the following table we list the pivotal player for all combinations of the input vector $a$ 
  and the ordering $\pi$. We assume that each player votes $l$ with probability $p_l$, i.e., $p_1,p_2\ge 0$ and $p_1+p_2=1$, and 
  that the votes of the two players are independent.
  \begin{center}
    \begin{tabular}{|rcccccc|r|}
      \hline
      $a\backslash \pi$ & $(1,2,3)$ & $(1,3,2)$ & $(2,1,3)$ & $(2,3,1)$ & $(3,1,2)$ & $(3,2,1)$ & $p(a)$\\
      \hline
      $(1,1,1)$ & 1 & 1 & 1 & 3 & 1 & 2 & $p_1^3$\\
      $(1,1,2)$ & 1 & 1 & 1 & 1 & 1 & 1 & $p_1^2p_2$\\
      $(1,2,1)$ & 1 & 1 & 1 & 1 & 1 & 1 & $p_1^2p_2$\\
      $(1,2,2)$ & 1 & 1 & 1 & 3 & 1 & 2 & $p_1p_2^2$\\
      $(2,1,1)$ & 3 & 2 & 3 & 3 & 2 & 2 & $p_1^2p_2$\\
      $(2,1,2)$ & 3 & 2 & 3 & 1 & 2 & 1 & $p_1p_2^2$\\
      $(2,2,1)$ & 3 & 2 & 3 & 1 & 2 & 1 & $p_1p_2^2$\\
      $(2,2,2)$ & 3 & 2 & 3 & 3 & 2 & 2 & $p_2^3$\\
      \hline
    \end{tabular}  
  \end{center}
  Given the probabilities the (generalized) influence measure for players $2$ and $3$ are given by $\frac{1+2p_2^2}{6}$ and 
  the (generalized) influence measure for player $1$ is given by $\frac{4-4p_2^2}{6}$, i.e., the values are not independent 
  from the probability distribution.
\end{Example} 

\section{The axiomatic approach}
\label{sec_axioms}

\noindent
The meaningfulness of influence measures or power indices is commonly justified by providing some axioms which are satisfied by 
the measure and uniquely determine it. Here we go along these lines for our influence measure from Definition~\ref{def_influence_measure}. 
To this end let $v$ be a non-trivial $(J,K)$ game on $n$ players with $|\operatorname{im}(v)|>1$. Due to Lemma~\ref{lemma_efficient} the 
influence measure $\varphi$ is \emph{efficient} for $v$. For any null player $i$ in $v$ we obviously have $\varphi_i(v)=0$, i.e., $\varphi$ 
satisfies the \emph{null player axiom}. For each $\pi\in S_n$ we define the game $\pi v$ by $(\pi v)(a)=v\!\left(a_{\pi(1)},\dots,a_{\pi(n)}\right)$. 
With this, we have $\varphi_i(\pi v)=\varphi_{\pi(i)}(v)$, which is called the \emph{anonymity axiom}. The \emph{transfer axiom} 
is satisfied if we have 
$$
  \varphi_i(u)+\varphi_i(w)= \varphi_i(u\vee w)+\varphi_i(u\wedge w)
$$    
for all $1\le i\le n$, where $(u\vee w)(a):=\max\{u(a),w(a)$ and $(u\wedge w)(a):=\min\{u(a),w(a)$. 

\begin{Lemma}
  Let $u,w$ be nontrivial, surjective, monotonic, output-rough $(J,K)$ games. Then, both 
  $u\vee w$ and $u\wedge w$ are nontrivial, surjective, monotonic, output-rough $(J,K)$ games.
\end{Lemma}
\begin{Proof}
  Obviously $u\vee w$ and $u\wedge w$ are nontrivial, surjective, monotonic $(J,K)$ games. 

  Next we prove that $u\wedge w$ is output-rough. To this end let $a=(a_1,\dots,a_n)\in J^n$, 
  $\pi\in S_n$, and $1\le h\le n$ arbitrary but fix. Choose integers $\alpha_1,\alpha_2,\beta_1,\beta_2$ 
  such that
  $$
    \left\{ u(a_1',\dots,a_n')\,:\, a_l'=a_l\text{ if } \pi(l)\le h \text{ and } a_l'\in J\text{ otherwise}\right\}=[\alpha_1,\beta_1]\cap\mathbb{Z}
  $$ 
  and
  $$
    \left\{ w(a_1',\dots,a_n')\,:\, a_l'=a_l\text{ if } \pi(l)\le h \text{ and } a_l'\in J\text{ otherwise}\right\}=[\alpha_2,\beta_2]\cap\mathbb{Z}.
  $$
  We will prove
  $$
    \left\{ (u\wedge w)(a_1',\dots,a_n')\,:\, a_l'=a_l\text{ if } \pi(l)\le h \text{ and } a_l'\in J\text{ otherwise}\right\}=
    \left[\alpha_3,\beta_3\right]\cap\mathbb{Z},
  $$
  where $\alpha_3=\min\{\alpha_1,\alpha_2\}$ and $\beta_3=\min\{\beta_1,\beta_2\}$, in the following. 
  We set 
  $$
    R=\left\{(a_1',\dots,a_n')\in J^n\,:\, a_l'=a_l\text{ if } \pi(l)\le h\right\}.
  $$  
  Now let $\gamma\in [\alpha_3,\beta_3]$. If $\gamma<\alpha_2$ we choose a vector $r\in R$ with $u(r)=\gamma$. Since $w(r)>\gamma$, we have 
  $(u\wedge w)(r)=\gamma$. Similarly, we can conclude the existence of an input vector $r\in R$ with $(u\wedge w)(r)=\gamma$ if 
  $\gamma<\alpha_1$, $\gamma>\beta_1$, or $\gamma>\beta_2$. In the remaining cases we have $\alpha_1\le\gamma\le\beta_1$ and 
  $\alpha_2\le\gamma\le\beta_2$. Now let $R_\gamma$ be the set of elements $r'\in R$ with $u(r')=\gamma$. If there exists an element 
  $r\in R_{\gamma}$ with $w(r)\ge\gamma$, then $(u\wedge w)(r)=\gamma$. So, we assume $w(r')<\gamma$ for all $r'\in R_\gamma$ and choose 
  an arbitrary $r^h\in R_\gamma$. For each $h<l\le n$ we define $r^l$ by $r^l_i=r^{l-1}_i$ if $\pi(i)!=l$ and 
  $r^l_i=j$ for $\pi(i)=l$. By construction we have $r^h<\dots<r^n=(j,\dots,j)$ and $w(r^n)\ge \gamma$. Let $g$ be the smallest index 
  with $w(r^g)\ge \gamma$ and $w(r^{g-1})<\gamma$. Since $w$ is output-rough we can modify the input of player $i$ with $\pi(i)=g$ 
  in $r^{g-1}$ to a vector $r$ such that $w(r)=\gamma$ and $r\ge r^{g-1}$. Since $u(r)\ge u(r^{g-1})\ge\gamma$ we have 
  $(u\wedge w)(r)=\gamma$.
  
  For $u\vee w$ we can similarly conclude
  $$
    \left\{ (u\vee w)(a_1',\dots,a_n')\,:\, a_l'=a_l\text{ if } \pi(l)\le h \text{ and } a_l'\in J\text{ otherwise}\right\}=
    \left[\alpha_4,\beta_4\right]\cap\mathbb{Z},
  $$
  where $\alpha_4=\max\{\alpha_1,\alpha_2\}$, $\beta_4=\max\{\beta_1,\beta_2\}$, and $\alpha_1,\alpha_2,\beta_1,\beta_2$ are chosen 
  as above.
\end{Proof}

From Lemma~\ref{lemma_influence_measure_monotonic} we conclude:
\begin{Lemma}
  For each non-trivial, surjective, monotonic, output-rough $(J,K)$ games $u$ and $w$ on $n$ players we have
  $\varphi_i(u)+\varphi_i(w)= \varphi_i(u\vee w)+\varphi_i(u\wedge w)$ for all $1\le i\le n$. 
\end{Lemma}

In other words, the influence measure $\varphi$ from Definition~\ref{def_influence_measure} satisfies the four 
classical axioms, used for the first axiomatization of the Shapley-Shubik index, see \cite{dubey1975uniqueness}, on 
the class of non-trivial, surjective, monotonic, output-rough $(J,K)$ games.

\begin{Lemma}
  Let $\gamma$ be a mapping from the set of non-trivial, surjective, monotonic $(j,2)$ games on $n$ players 
  to $\Delta_n:=\left\{(x_1,\dots,x_n)\in [0,1]^n\,:\, \sum_{i=1}^n x_i=1\right\}$. If $\gamma$ satisfies the
  transfer axiom, then $\gamma(v)$ can be recursively computed from the values of $\gamma$ for games with a 
  unique minimal winning coalition.
\end{Lemma}
\begin{Proof}
  We prove by induction on the number of minimal winning coalitions. By $u_a$ we denote the game with unique 
  minimal winning $a$. If the (pairwise different) minimal winning vectors of $v$ are given by $a^1,\dots,a^l$, 
  then we can write $v=x\vee y$, where $x=u_{a^1}\vee \dots u_{a^{l-1}}$ and $y=u_{a^l}$, so that $x$ has $l-1$ 
  and $y$ has $1$ minimal winning vector. Since $x\wedge y$ has at most $l-1$ winning vectors, we can compute recursively 
  compute $\gamma(v)=\gamma_(x)+\gamma(y)-\gamma(x\wedge y)$.
\end{Proof}

For non-trivial, surjective, monotonic $(2,2)$ games anonymity, efficiency, and the null player axiom uniquely determine 
the value of $\gamma$ on each game consisting of a single minimal winning vector. Here the null players obtain $\gamma_i=0$ 
and the non-null players obtain one divided by the number of non-null players. For $(j,2)$ games with $j>2$ the situation 
is more involved for our influence measure from Definition~\ref{def_influence_measure}:

\begin{Example}
  Let $v$ the the surjective, monotonic $(4,2)$ game on $3$ players with unique minimal winning vector $(2,3,4)$. Next we 
  determine the number of cases $(a,\pi)$, where each player is pivotal. For player~$1$ we consider the cases:
  \begin{itemize}
    \item $\pi=(1,\star,\star)$: $a_1=1$ $\rightarrow$ $2\cdot 1\cdot 4\cdot 4=32$ cases;
    \item $\pi=(2,1,3)$: $a_1=1$, $a_2\in\{3,4\}$ $\rightarrow$ $1\cdot 1\cdot 2\cdot 4=8$ cases;
    \item $\pi=(3,1,2)$: $a_1=1$, $a_3=4$ $\rightarrow$ $1\cdot 1\cdot 4\cdot 1=4$ cases;
    \item $\pi=(\star,\star,1)$: $a_2\in\{3,4\}$, $a_3=4$ $\rightarrow$ $2\cdot 4\cdot 2\cdot 1=16$ cases.
  \end{itemize} 
  For player~$2$ we consider the cases:
  \begin{itemize}
    \item $\pi=(2,\star,\star)$: $a_2\in\{1,2\}$ $\rightarrow$ $2\cdot 4\cdot 2\cdot 4=64$ cases;
    \item $\pi=(1,2,3)$: $a_1=1$, $a_1\in\{2,3,4\}$, $a_2\in\{1,2\}$ $\rightarrow$ $1\cdot 3\cdot 2\cdot 4=24$ cases;
    \item $\pi=(3,2,1)$: $a_1=1$, $a_2\in\{1,2\}$, $a_3=4$ $\rightarrow$ $1\cdot 4\cdot 2\cdot 1=8$ cases;
    \item $\pi=(\star,\star,2)$: $a_1\in\{2,3,4\}$, $a_3=4$ $\rightarrow$ $2\cdot 3\cdot 4\cdot 1=24$ cases.
  \end{itemize}
  For player~$3$ we consider the cases:
  \begin{itemize}
    \item $\pi=(3,\star,\star)$: $a_3\in\{1,2,4\}$ $\rightarrow$ $2\cdot 4\cdot 4\cdot 3=96$ cases;
    \item $\pi=(1,3,2)$: $a_1=1$, $a_1\in\{2,3,4\}$, $a_3\in\{1,2,3\}$ $\rightarrow$ $1\cdot 3\cdot 4\cdot 3=36$ cases;
    \item $\pi=(2,3,1)$: $a_1=1$, $a_2\in\{3,4\}$, $a_3\in\{1,2,3\}$ $\rightarrow$ $1\cdot 4\cdot 2\cdot 3=24$ cases;
    \item $\pi=(\star,\star,3)$: $a_1\in\{2,3,4\}$, $a_2\in\{3,4\}$ $\rightarrow$ $2\cdot 3\cdot 2\cdot 4=48$ cases.
  \end{itemize}
  Thus, we have $\varphi_{1}(v)=\frac{60}{386}=\frac{5}{32}$, $\varphi_{2}(v)=\frac{120}{386}=\frac{10}{32}=\frac{5}{16}$, 
  and $\varphi_{3}(v)=\frac{204}{386}=\frac{17}{32}$. 
\end{Example} 

An additional property of our influence measure $\varphi$ is that for a game $v'$ arising from $v$ by adding a null player $l$, 
we have $\varphi_i(v')=\varphi_i(v)$ for all $i\neq l$ (and $\varphi_l(v')=0$).

\section{The limiting case}
\label{sec_limiting_case}

\noindent
By the following normalization trick we can remove the $\frac{1}{k-1}$ factor in Equation~(\ref{eq_influence_measure_monotonic}) 
in Lemma~\ref{lemma_influence_measure_monotonic}. Instead of $K=\{1,\dots,k\}$ we use 
$K=\left\{\frac{0}{k-1},\frac{1}{k-1},\dots,\frac{k-1}{k-1}\right\}\subseteq[0,1]$. Then the $\underline{v}$- and $\overline{v}$-values 
contain the necessary factor itself.\footnote{By choosing another segmentation one can implement the evaluation function 
proposed in \cite{freixas2005shapley} to assign different \textit{weights} to the output values.} 
Without any substantial effect we may also relabel the set of inputs from $J=\{1,\dots, j\}$ to
$J=\left\{\frac{0}{j-1},\frac{1}{j-1},\dots,\frac{j-1}{j-1}\right\}\subseteq[0,1]$. 
However, by simultaneously increasing 
$j$ and $k$ (possibly at different velocities) we obtain an approximation of a voting scheme $[0,1]^n\rightarrow [0,1]$. 
We now introduce those objects directly.

\begin{Definition}
 A \emph{$([0,1],[0,1])$ game} on $n$ players is a mapping $v:[0,1]^n\rightarrow [0,1]$. We call $v$ \emph{surjective}, 
 \emph{monotonic}, or \emph{continuous} if the mapping is surjective, weakly monotonic increasing, or continuous, respectively. 
\end{Definition} 

As an abbreviation, we speak of a \emph{continuous game} $v$ if $v$ is a surjective, monotonic, continuous $([0,1],[0,1])$ game. 
Similarly to Definition~\ref{def_null_player} we call a player $i$ \emph{null player} in $v$, if $v(x_1,\dots,x_n)=
v(x_1,\dots,x_{i-1},x_i',x_{i+1},\dots,x_n)$ for all $x_1,\dots,x_n,x_i'\in[0,1]$.

As the continuity in a continuous game plays the role of output-roughness in a $(j,k)$ game, we reformulate 
Lemma~\ref{lemma_influence_measure_monotonic} to:
\begin{Definition}
  \label{def_influence_measure_continous}
  For each continuous game $v$ on $n$ players we define $\varphi_i(v)$ by
  $$
    \frac{1}{n!}\sum_{\pi\in S_n}\!\!\int_0^1\!\!\dots\!\!\int_0^1\!\!
    \left(\overline{v}_{\pi,\pi^{-1}(i)\!-\!1}(a)\!-\!\underline{v}_{\pi,\pi^{-1}(i)\!-\!1}(a)\right)
     -\left(\overline{v}_{\pi,\pi^{-1}(i)}(a)\!-\!\underline{v}_{\pi,\pi^{-1}(i)}(a)\right) 
     \operatorname{d}x_1\dots \operatorname{d}x_n, 
  $$
  where $\overline{v}_{\pi,h}(a)=v(a_1',\dots,a_n')$ with $a_l'=a_l$ if $\pi(l)\le h$ and $a_l'=1$ otherwise;
  $\underline{v}_h(b)=v(b_1',\dots,b_n')$ with $b_l'=b_l$ if $\pi(l)\le h$ and $b_l'=0$ otherwise, for all $1\le i\le n$.
\end{Definition}

We remark that the influence measure from Definition~\ref{def_influence_measure_continous} satisfies efficiency, anonymity, 
the null player and the transfer axiom. Given a continuous game $v$ one can construct a series of monotonic 
$(j,j)$ games $v^{j}$, where $j\to\infty$, such that the values of $\varphi_i(v^{j})$ tend to $\varphi_i(v)$ for all $1\le i\le n$.

For the examples $\hat{v}(x_1,x_2,x_3)=\frac{1x_1^2+2x_2^2+3x_3^2}{6}$ and $\tilde{v}(x_1,x_2,x_3)=x_1x_2^2x_3^3$ the function 
$\varphi$ was evaluated in \cite{kurz2014measuring}: $$\varphi(\hat{v})=\left(\frac{1}{6},\frac{2}{6},\frac{3}{6}\right)=
(0.1\overline{6},0.\overline{3},0.5)$$ and $$\varphi(\tilde{v})=\left(\frac{35}{144},\frac{50}{144},\frac{59}{144}\right)=
(0.2430\overline{5},0.347\overline{2},0.4097\overline{2}).$$  

While typically the evaluation of $\varphi$ for a continuous game is based on rather tedious case distinctions, there are simple 
formulas for special cases, see \cite{kurz2018influence} for more.

%% Stimmt das unten stehende Theorem noch fuer mehr Beispiele?
%% Falls ja, ist das definitiv einen komplizierteren Beweis wert und wuerde einen schoenen Bogen spannen
%% 
%% \begin{Definition}
%%   \label{def_weighted_median}
%%   Let $n$ be a positive integer and $w_1,\dots,w_n\in\mathbb{R}_{\ge 0}$ such that there exists no $S\subseteq N$ 
%%   with $w(S)=\frac{1}{2}\cdot w(N)$. By $v^{[w(N)/2;w_1,\dots,w_n]}$ we denote the continuous game that maps $x\in [0,1]^n$ 
%%   to the weighted median according to the weights $w_1,\dots, w_n$.
%% \end{Definition}
%% 
%% As an example we consider $v^{[5.5;5,3,2,1]}(0.7,0.2,0.5,0.1)=0.5$.
%% 
%% \begin{Theorem}
%%   Let $\tilde{v}=[q;w_1,\dots,w_n]$ be a simple game with $q=w(N)$ satisfying the conditions from Definition~\ref{def_weighted_median}, 
%%   defined by $\tilde{v}(S)=1$ iff $w(S)\ge q$ and $w(S)=0$ otherwise, i.e., $\tilde{v}$ is a so-called weighted game. For each $1\le i\le n$ 
%%   we have
%%   $$
%%     \tilde{\varphi}_i(\tilde{v})=\varphi_i\!\left(v^{[q;w_1,\dots,w_n]}\right).
%%   $$ 
%% \end{Theorem}
%% \begin{Proof}
%%   For arbitrary but fixed $x\in[0,1]^n$ with $x_j\neq x_h$ for $j\neq h$, let $i$ be the unique player whose input $x_i$ attains 
%%   the weighted median. With this we have $v^{[q;w_1,\dots,w_n]}(x)=x_i$.
%%   \bigskip
%%   \bigskip  
%%   \begin{center}
%%     \dots
%%   \end{center}
%%   \bigskip
%%   \bigskip
%%   $$
%%     \int_0^1 x^a(1-x)^{n-a-1}=\frac{a!(n-a-1)!}{n!},
%%   $$  
%%   see the well studied Euler beta function.
%% \end{Proof}

\begin{Theorem}
  Let $w_1,\dots,w_n\in[0,1]$ with $\sum_{i=1}^n w_i=1$ and $f_i:[0,1]\rightarrow[0,1]$ continuous monotonic functions 
  with $f_i(0)=0$ and $f_i(1)=1$ for all $1\le i\le n$. Then, $v(x)=\sum_{i=1}^n w_i\cdot f_i(x_i)$ defines a 
  continuous game and we have $\varphi_i(v)=w_i$ for all $1\le i\le n$.
\end{Theorem}
\begin{Proof}
  Let $q\le i\le n$ and $\pi\in S_n$ be arbitrary but fixed. With $S=\{l\in N\,:\, \pi(l)<\pi(i)\}$ we have
  \begin{eqnarray*}
    \overline{v}_{\pi,\pi^{-1}(i)-1}(x) &=& \sum_{l\in S} w_l\cdot f_l(x_l)\,+\, \sum_{l\in N\backslash S} w_l\cdot f_l(1)
     =  \sum_{l\in S} w_l\cdot f_l(x_l) \,+\, 1-w(S)\\
    \overline{v}_{\pi,\pi^{-1}(i)}(x) &=& \sum_{l\in S} w_l\cdot f_l(x_l)\,+\,w_i\cdot f_i(x_i)\,+\, 1-w(S)-w_i\\
    \underline{v}_{\pi,\pi^{-1}(i)-1}(x) &=& \sum_{l\in S} w_l\cdot f_l(x_l)\,+\, \sum_{l\in N\backslash S} w_l\cdot f_l(0)
     =  \sum_{l\in S} w_l\cdot f_l(x_l)\\
     \underline{v}_{\pi,\pi^{-1}(i)}(x) &=& \sum_{l\in S} w_l\cdot f_l(x_l)\,+\,w_i\cdot f_i(x_i)\\
  \end{eqnarray*}
  so that
  $
    \left(\overline{v}_{\pi,\pi^{-1}(i)\!-\!1}(a)\!-\!\underline{v}_{\pi,\pi^{-1}(i)\!-\!1}(a)\right)
     -\left(\overline{v}_{\pi,\pi^{-1}(i)}(a)\!-\!\underline{v}_{\pi,\pi^{-1}(i)}(a)\right)
     =w_i
  $.
\end{Proof}

So the influence distribution for the example $\hat{v}$ is no surprise, while for continuous games similar to $\tilde{v}$ 
no general and easy to evaluate formula is known.

The analogy to different probabilities for the discrete set of input states is a density function in the continuous case.
An example with density functions $f_1(x)=\frac{3}{4}\cdot(1-x^2)$, $f_2(x)=f_2(x)=\frac{3}{8}\cdot(1+x^2)$ and $v(x)$ equal 
to the median of $x_1,x_2,x_3$ was computed in \cite{kurz2014measuring}: $\varphi_1(v)=\frac{554}{13440}\approx 0.04122$ and 
$\varphi_2(v)=\varphi_3(v)=\frac{563}{13440}\approx 0.04189$.

\section{Conclusion}
\label{sec_conclusion}
      
\noindent
We have studied the classical Shapley-Shubik index for simple games in the roll-call model from \cite{felsenthal1996alternative} 
and gave a direct combinatorial proof for the fact that the expected number of cases where a certain player is pivotal is independent 
from the specific distribution of the {\lq\lq}yes{\rq\rq}- and {\lq\lq}no{\rq\rq}-votes, as long as the probability 
does only depended on the number of {\lq\lq}yes{\rq\rq}-votes. This generalizes the result from \cite{felsenthal1996alternative} and 
gives a nice and vivid description for the Shapley-Shubik index that is less artificial than previous ones. We have applied the 
roll-call model for the generalized case of $j\ge 2$ ordered input and $k\ge 2$ ordered output states. Some notation from simple 
games can be generalized in a consistent way. By considering the reduction of uncertainty we have tried to provide a more 
persuasive basis for the generalized Shapley-Shubik index from \cite{freixas2005shapley} and other places. We do not claim 
that this influence measure is the correct generalization of the Shapley-Shubik index and we agree with the authors of 
\cite{felsenthal2013models} that the status of the Shapley-Shubik index for ternary and more general games requires further study. 
For further generalizations of the Shapley-Shubik index we refer the interested reader to \cite{roth1988shapley}.  
A first set of results with respect to the axiomatic approach is presented in order justify the proposed influence measure by another 
approach. In any case we find it advantageous to have an influence measure at hand that can, to some extend, be consistently defined 
for larger classes of games with several levels of approval in the input and output. Even continuous input and output spaces 
make sense in decisions on e.g.\ rate of taxes or other continuous variables.      
 
%% \begin{acknowledgement}
%%   An acknowledgement may be placed at the end of the article.
%% \end{acknowledgement}

%%\bibliographystyle{amsplain}
%%\bibliographystyle{plain}
%%\bibliography{Roll-call}

\end{document}